\documentclass[a4paper,12pt]{article}

\usepackage{color}

\usepackage{bm}
 \usepackage{amssymb}
 \usepackage{amsbsy}
 \usepackage{amscd}
 \usepackage{amsmath}
 \usepackage{amsthm}
\usepackage{amsmath}
\usepackage{amsfonts}

\newcommand{\RR}{\mathbb R}
\newcommand{\CC}{\mathbb C}
\newcommand{\EE}{\mathbb E}
\newcommand{\VV}{\mathbb V}

\begin{document}

\title{
Recent Advances in Algebraic Geometry and Bayesian Statistics
}

\author{Sumio Watanabe
\\
Department of Mathematical and Computing Science
\\
Tokyo Institute of Technology 
\\
2-12-1 Oookayama, Meguro-ku,Tokyo  52-8552 Japan. \\
mailbox W8-42\\
E-mail: swatanab@c.titech.ac.jp
\\
ORCID:0000-0001-8341-5639
}

\date{}

\maketitle

\begin{abstract}

This article is a review of theoretical advances in the research field of 
algebraic geometry and Bayesian statistics in the last two decades. 

Many statistical models and learning machines which contain
hierarchical structures or latent variables are 
called nonidentifiable, because 
the map from a parameter to a statistical model is not
one-to-one. In nonidentifiable models, 
both the likelihood function and the posterior distribution 
have singularities in general, 
hence it was difficult to analyze their statistical properties. 
However, from the end of the 20th century, 
new theory and methodology based on algebraic geometry 
 have been established which enables 
us to investigate such models and machines in the real world. 

In this article, the following results in recent advances are reported. 
First, we explain the framework of Bayesian statistics and
introduce a new perspective from the birational geometry. 
Second, two mathematical solutions are derived 
based on algebraic geometry. 
An appropriate parameter space can be found by a resolution map, 
which makes the posterior distribution be normal crossing and 
the log likelihood ratio function be well-defined. 
Third, three applications to statistics are introduced. 
The posterior distribution is represented by the renormalized form, 
the asymptotic free energy is derived, and the universal
formula among the generalization loss, the cross validation, and
the information criterion is established. 

Two mathematical solutions and three applications to statistics
based on algebraic geometry reported in this article are now being
used in many practical fields in data science and artificial intelligence.
\end{abstract}

\section{Introduction}

Many statistical models and learning machines which have 
hierarchical structure or latent variables are widely used 
 in data science, artificial intelligence, bioinformatics, economics, 
political science, psychology, and so on. 
Such models and machines are called {\it nonidentifiable}, 
since the map 
from a parameter to a probability density function is not one-to-one 
\cite{Watanabe2001a}. 
They are also called {\it singular}, because
both the likelihood function and the posterior distribution
contain singularities which cannot 
be approximated by any Gaussian function.
In fact, the Fisher information matrix contains zero eigen values
 and the Laplace approximation 
does not capture the essential property of the posterior distribution. 
In other words, 
the classical statistical theory which needs the regularity condition 
does not hold, resulting that 
hypothesis test, model selection, hyperparameter optimization, and learning by 
gradient decent should be studied from new perspective 
 \cite{Hartigan1985,Hagiwara1993,Watanabe1995,Fukumizu1996,Amari2001}. 

Such nonidentifiable and singular models and machines are not special but 
ubiquitous in modern statistics and machine learning. 
For example, neural networks and deep learning 
\cite{Watanabe2001,Watanabe2001a,Aoyagi2012}
have hierarchical structures, 
normal mixtures \cite{Yamazaki2003,Kariya2020},
Poisson mixtures \cite{Sato2019}, 
multinomial mixtures \cite{WatanabeT2022}, 
and latent Dirichlet allocations \cite{Hayashi2021} 
have latent or hidden variables, and 
matrix factorizations \cite{Aoyagi2005,Hayashi2017}, 
Boltzmann machines \cite{Yamazaki2005a,Aoyagi2010}, and
Markov models \cite{Yamazaki2005,Zwie2011} have 
both hidden and hierarchical parts. 
In other words, almost all statistical models and 
learning machines which extract hidden structures or hierarchical inferences
are nonidentifiable and singular \cite{Watanabe2007}. 

Here let us illustrate the two 
mathematical problems which are universally found in
nonidentifiable and singular models and machines. 
Let $X$ and $X_1,X_2,...,X_n$ be independently and identically distributed 
$\RR^N$-valued random variables and let 
$h(x,\theta)$ be a real-valued function on $\RR^N\times\RR^d$ which 
is analytic for $\theta$ and satisfies $H(\theta)\equiv\EE[h(X,\theta)]\geq 0$.
In statistics and machine learning, 
a random process  $H_n(\theta)$ defined on $\RR^d$, 
\[
H_n(\theta)=\frac{1}{n}\sum_{i=1}^n h(X_i,\theta)
\]
is quite often studied, in fact, the 
minus log density ratio function is a typical example. 
It can be rewritten as 
\[
nH_n(\theta)=nH(\theta)-(nH(\theta))^{1/2}A_n(\theta),
\]
where $A_n(\theta)$ is defined by 
\[
A_n(\theta)=\frac{1}{\sqrt{n}} \sum_{i=1}^n 
\left\{\frac{H(\theta)-h(X_i,\theta)}{H(\theta)^{1/2}}\right\}.
\]
Note that $\EE[A_n(\theta)]=0$ for any $\theta$ by the definition. 
If $H(\theta)=0$ at a unique parameter $\theta_0$ and the Hessian matrix
of $H(\theta)$ at $\theta_0$ is positive definite, then regular statistical theory 
holds,  resulting that the learning curves are determined by the
 dimension of the parameter 
\cite{Amari1992,Amari1993,Murata1995}. 
However, if otherwise, there are two mathematical difficulties. 

The first mathematical problem in singular statistics 
is that the set defined by $H(\theta)=0$ consists of not a single element 
but multiple and uncountable elements with singularities. 
Hence the function $H(\theta)$ cannot be 
approximated by any quadratic form in neighborhoods of $H(\theta)=0$, 
moreover, it is difficult to measure 
the volume of 
$\{\theta\;;\;H(\theta)<t\}$ for $t\rightarrow +0$ near singularities
although it determines the accuracy of Bayesian statistics. 

The second problem is that the random process $A_n(\theta)$ is 
not well-defined on $H(\theta)=0$. In the region
 $H(\theta)>0$, $A_n(\theta)$ is an empirical process which
 converges to a Guassian process in distribution 
as $n$ tends to infinity. However, $H(\theta)-h(x,\theta)$ is not
divisible by $H(\theta)^{1/2}$ in general, hence 
in a neighborhood of $H(\theta)\rightarrow +0$, the limit process 
$A_n(\theta)\rightarrow 0/0$ cannot be determined uniquely. 

In this paper, we show that algebraic geometry
gives a mathematically natural solution to the foregoing
statistically essential problems. 
Assume that $H(\theta)$ is an arbitrary analytic function. Then, 
based on the resolution theorem which was proved by Hironaka 
\cite{Hironaka1964} and applied by Atiyah and Kashiwara \cite{Atiyah1970,Kashiwara1976}, 
there exists an appropriate analytic function from a manifold to the parameter space, 
\[
\theta=g(u), 
\]
such that $H(g(u))$ is {\it normal crossing} on an arbitrary neighborhood of
$H(g(u))=0$. Here a function $H(g(u))$ is said to be
normal crossing if $H(g(u))$ is represented by a simple direct 
product of $(u_j)^{k_j}$, 
whose mathematical definition is given  in Section \ref{section:alg}. 
This theorem enables us to study the function $H_n(\theta)$ and solves two mathematical
problems. First, we can derive that 
the volume of the set 
$
\{w;H(w)<t\} $ as $t\rightarrow +0$ 
is in proportion to 
$t^{\lambda}(-\log t)^{m-1}$, where $\lambda$ and $m$ are birational
invariants called a real log canonical threshold and a multiplicity. 
We show that the accuracy of Bayesian statistics is determined by this
volume. 
Second, also we show that 
 $A_n(g(u))$ is made to be a well-defined function 
of $u$ because $H(g(u))-h(x,g(u))$ is divisible by $H(g(u))^{1/2}$,
 hence the empirical process theory ensures
that $A_n(g(u))$ converges to a Gaussian process on $u$ in distribution
even in an arbitrary neighborhood of $H(g(u))=0$. 

These two solutions have three applications to statistics. First, 
the posterior distribution is represented by the renormalized form, which 
clarifies the scaling law of the posterior distribution. 
Second, 
the asymptotic free energy, which is equal to
the minus log marginal likelihood, is derived.  Lastly,  the universal
formula among the generalization loss, the cross validation, and
the information criterion is proved, by which the generalization loss for 
unknown information source can be estimated. 

This article consists of six chapters.
In Section \ref{section:basic}, we explain the framework of 
Bayesian statistics and introduce two important random variables, 
the free energy and the generalization loss, which are measures of
appropriateness of a statistical model and a prior distribution. 
In Section \ref{section:alg}, mathematical description of Hironaka resolution
theorem is introduced, by which we can represent the average 
log density ratio function can be made normal crossing, and 
two birational invariants, the real log canonical 
threshold and its multiplicity are defined.
In Sections \ref{section:two} and \ref{section:three}, 
two mathematical solutions and three applications to statistics are 
explained. 
Lastly, in Section \ref{section:conc}, the results of this paper is summarized.

\section{Framework of Bayesian Statistics}\label{section:basic}

In this section, we prepare a mathematical framework 
of Bayesian statistics and explain the purpose of this paper. 

Let $p(x|\theta)$ be a conditional probability density function 
of $x\in\RR^N$ for
a given parameter $\theta\in\Theta\subset\RR^d$
and $\pi(\theta)$ be a probability density function of $\theta$. 
In Bayesian statistics, a candidate pair made by a person
\begin{align}
\theta&\sim \pi(\theta), \label{eq:prior}
\\
X^n&\sim \prod_{i=1}^n p(x_i|\theta),\label{eq:model}
\end{align}
is investigated, which means
$X^n=\{X_1,X_2,...,X_n\}$ may be subject to a statistical model $p(x|\theta)$
with a prior distribution $\pi(\theta)$. 

If $X_1,X_2,...,X_n,...,$ are assumed to be subject to some
distributions such as eqs. (\ref{eq:prior}) and (\ref{eq:model}), 
they are called exchangeable. If they are exchangeable, 
 by de Finetti's theorem, there exist both
 ${\cal Q}(q)$ and $q(x)$ such that
\begin{align}
q(x)&\sim {\cal Q}(q),\label{eq:prior0}
\\
X^n&\sim \prod_{i=1}^n q(x_i),\label{eq:model0}
\end{align}
where ${\cal Q}(q)$ is a probability distribution on
the set of all probability distributions on $\RR^N$ and $q(x)$ is
 a probability density function which is subject to ${\cal Q}(q)$. 
The general pair ${\cal Q}(q)$ and $q(x)$ contains
a specific pair $\pi(\theta)$ and $p(x|\theta)$, hence
if a person makes a candidate pair $\pi(\theta)$ and $p(x|\theta)$ 
and rejects the existence of unknown ${\cal Q}(q)$ and $q(x)$, 
it is a mathematical contradiction \cite{Watanabe2022, Watanabe2022a}. 

In statistical science of a large world, both ${\cal Q}(q)$ and $q(x)$ are unknown
and all models are wrong \cite{Box1976}. A person cannot believe in
a specific pair of a statistical model and 
a prior distribution because it is under- or  over-parametrized in 
an environment of unknown uncertainty \cite{GelmanBDA,Gelman2013,Binmore2017,McElreath2020}. 
In other words, a person who made $p(x|\theta)$ and $\pi(\theta)$ is 
aware that both are only fictional candidates, 
resulting that it is necessary to check or evaluate its appropriateness 
from a mathematically general viewpoint \cite{Akaike1974,Akaike1980,Watanabe2022a}. 
In an older Bayesian statistics, a person needs to believe that 
${\cal Q}(q)=\pi(\theta)$ and $q(x)=p(x|\theta)$, whereas 
in modern Bayesian statistics, a person is aware to distinguish the 
unknown data-generating process
from a model and a prior \cite{Watanabe2022a}.

In this paper, we explain a role of algebraic geometry in modern Bayesian statistics. 
It is assumed that 
$X^n$ and $X$ are independent and generated from the unknown data-generating 
process $q(x)$ which is subject to ${\cal Q}(q)$ and that 
a candidate pair, $p(x|\theta)$ and $\pi(\theta)$, is prepared by a person. 
Let $\EE[f(X^n)|q]$ and $\EE_X[f(X)|q]$ denote the expectation values of  given
 functions $f(X^n)$ and $f(X)$ according to $\prod_i q(x_i)$ and $q(x)$, respectively. 

The average and empirical log loss
functions $L(\theta)$ and $L_n(\theta)$ are defined by
\begin{align}
L(\theta)&=-\EE_X[\log p(X|\theta)|q], 
\label{eq:L}
\\
L_n(\theta)&=-\frac{1}{n}\sum_{i=1}^n \log p(X_i|\theta). 
\label{eq:Ln}
\end{align}
The posterior distribution and the posterior 
predictive distribution using a candidate pair
eqs. (\ref{eq:prior}) and (\ref{eq:model}) are defined
respectively by
\begin{align}
p(\theta|X^n)&=\frac{1}{p(X^n)}\pi(\theta)\prod_{i=1}^n p(X_i|\theta),
\\
p(x|X^n)&=\int p(x|\theta)p(\theta|X^n)d \theta,
\end{align}
where 
\[
p(X^n)=\int\pi(\theta)\prod_{i=1}^n p(X_i|\theta)d\theta
\]
is the marginal likelihood. The average and variance by using 
the posterior distribution $p(\theta|X^n)$ are denoted by $\EE_\theta[\;\;]$ and
$\VV_\theta[\;\;]$ respectively. 

The free energy $F_n$, which is
equal to the minus log 
marginal likelihood, and the generalization loss $G_n$ are
defined respectively by
\begin{align}
F_n&=-\log p(X^n),\label{eq:Fn}
\\
G_n&=-\EE_X[ \log p(X|X^n)|q].\label{eq:Gn}
\end{align}
Then it follows that
\begin{align}
\EE[F_n|q]&={\rm KL}(q(X^n)||p(X^n))+nS(q),
\\
\EE[G_n|q]&={\rm KL}(q(X)||p(X|X^n))+S(q),
\end{align}
where ${\rm KL}(\;\;||\;\;)$ is Kullback-Leibler divergence 
of $q_1(x)$ and $q_2(x)$, 
\[
{\rm KL}(q_1(X)||q_2(X))=\int q_1(x)\log\frac{q_1(x)}{q_2(x)}dx,
\]
and $S(q)$ is the entropy of $q(x)$, 
\[
S(q)=-\int q(x)\log q(x)dx.
\]
Note that $S(q)$ does not depend on the candidate pair $p(x|\theta)$ and
$\pi(\theta)$. 
Therefore, the average free energy and the average 
generalization loss are minimized if and only if
$q(X^n)=p(X^n)$ and $q(x)=p(x|X^n)$, respectively. These properties
show that the free energy and the generalization loss 
can be understood as different measures of appropriateness of 
the pair $p(x|\theta)$ and $\pi(\theta)$. 
By the definition,  for an arbitrary positive
integer $n$,
\[
\EE[G_n|q]=\EE[F_{n+1}|q]-\EE[F_n|q]
\]
holds, however, the pair $(p(x|\theta),\pi(\theta))$ that
minimizes $F_n$ is different from the pair that minimizes
$G_n$ \cite{Watanabe2022a}. 

Minimizing the free energy, which is equivalent to 
maximizing the marginal likelihood, is often employed in 
Bayesian model selection and hyperparameter optimization \cite{Akaike1980}. 
It is also known that the difference of the free energies 
between the null hypothesis pair and 
the alternative pair gives the most powerful Bayesian test \cite{Kariya2020}. 
Also minimizing the generalization loss is often adopted for the purpose of 
the accurate prediction in statistics and machine learning 
\cite{Akaike1974,Watanabe2010}. Hence it is important to clarify the
mathematical properties of the free energy and the generalization loss. 

In order to estimate the generalization loss, three
random variables are defined, 
the training loss $T_n$, the leave-one-out cross validation $C_n$ \cite{Gelfand1992,Vehtari2002,Gelman2014}, and 
the widely applicable information criterion $W_n$ \cite{Watanabe2010} respectively by 
\begin{align}
T_n&=-\frac{1}{n}\sum_{i=1}^n\log p(X_i|X^n),\label{eq:Tn}
\\
C_n&=-\frac{1}{n}\sum_{i=1}^n\log p(X_i|X^n\setminus X_i),\label{eq:Cn}
\\
W_n&=T_n+\frac{1}{n}\sum_{i=1}^n \VV_\theta[\log p(X_i|\theta)],\label{eq:Wn}
\end{align}
where $X^n\setminus X_i$ is the set leaving $X_i$ out from $X^n$. 
We define a function $\Omega(\theta)$ on $\Theta$ by 
\begin{align}
\Omega(\theta) &\equiv\pi(\theta)\prod_{i=1}^n p(X_i|\theta).
\end{align}
Then the posterior distribution is represented by
\[
p(\theta|X^n) =\frac{\Omega(\theta)}
{\int \Omega(\theta')d\theta'},
\]
and the free energy is given by
\begin{align}
F_n&=-\log \int \Omega(\theta) d\theta. \label{eq:FOmega}
\end{align}
The posterior average $\EE_\theta[\;\;]$ and
the posterior variance $\VV_\theta[\;\;]$ are also represented by $\Omega(\theta)$, 
hence
\begin{align}
G_n&=-\EE_X[ \log \EE_\theta[p(X|\theta)]|q],\label{eq:GOmega}
\\
T_n&=-\frac{1}{n}\sum_{i=1}^n\log \EE_\theta[p(X_i|\theta)],\label{eq:Tn2}
\\
C_n&=\frac{1}{n}\sum_{i=1}^n\log \EE_\theta[1/p(X_i|\theta)],\label{eq:Cn2}
\\
W_n&=T_n+\frac{1}{n}\sum_{i=1}^n\log \VV_\theta[\log p(X_i|\theta)],\label{eq:Wn2}
\end{align}
are all represented by using $\Omega(\theta)$. 
If a statistical model contains hierarchical structure or latent variables, 
then $\Omega(\theta)$ cannot be 
approximated by any Gaussian function in general, because
a statistical model may be under- or over-parametrized. 
The main purpose of this paper is to clarify what mathematical structure determines 
$\Omega(\theta)$, 
the posterior distribution, and the five random variables.
\vskip3mm\noindent
{\bf Purpose of this paper.}
In this paper, we characterize the
integration $\Omega(\theta)d\theta$ from the algebro-geometric point of view,
and derive the probabilistic behaviors of five random variables, 
the free energy $F_n$, the generalization loss $G_n$,
the training loss $T_n$, the leave-one-out cross validation $C_n$, 
and the widely applicable information criterion $W_n$, 
even if $\Omega(\theta)$ is far from any Gaussian function. 
\vskip3mm\noindent
In this paper, we mainly study 
 the accuracy of the free energy and the posterior
predictive distribution.
Their performance for the case when the probabilities distributions are different 
between training and test are clarified \cite{Yamazaki2007} and 
the accuracy of the estimation
of the latent variables are also derived by \cite{Yamazaki2016}. 
In Bayesian statistics, it is one of the most important researches
 how to approximate the posterior 
distribution using Markov chain Monte Carlo (MCMC). Although the purpose 
of this paper is not studying MCMC methods, an 
algebro-geometric study of the posterior distribution may be useful in 
the design of MCMC process in singular models and machines. 
For example, the optimal sequence of the
 inverse temperatures in exchange Monte Carlo 
are clarified by the property of singular posterior distributions \cite{Nagata2008}.

\section{Algebro-Geometric Foundation}\label{section:alg}

In this section, we introduce an algebro-geometric foundation on 
which mathematical and statistical theory is constructed. 

For simplicity, we assume $\Theta\subset \RR^d$ is a compact set whose
open kernel is not empty 
and $L(\theta)$ is an analytic function of $\theta$ in some open set 
which contains $\Theta$. 
The set of all parameters that make $L(\theta)$ minimum is
\[
\Theta_0=\{\theta\in\Theta\;;\;L(\theta)\mbox{ is minimum}\}.
\]
If there exists $\theta_0\in\Theta_0$ such that
$q(x)=p(x|\theta_0)$, then $q(x)$ is said to be realizable by $p(x|\theta)$, 
or if otherwise it is said unrealizable. 
If $\Theta_0$ consists of a single element $\theta_0$ and if 
the Hessian matrix $\nabla^2 L(\theta_0)$ is positive definite, 
then $q(x)$ is said to be regular for $p(x|\theta)$, or 
if otherwise it is said to be singular. 

The set $\Theta_0$ is called an analytic set because it is the set
of all zero points of an analytic function $L(\theta)-L(\theta_0)$. If $L(\theta)$ is 
a polynomial function, then it is called an algebraic set. 
We assume that $p(x|\theta_0)$ does not depend on the choice of 
$\theta_0 \in\Theta$. If $q(x)$ is realizable by or regular for $p(x|\theta)$,
then such a condition is satisfied \cite{Watanabe2018}. For the case when 
$p(x|\theta_0)$ depends on $\theta_0\in\Theta_0$, see 
\cite{Nagayasu2022,Watanabe2010a}.

A log density ratio function $f(x,\theta)$ is 
defined by
\begin{align}\label{eq:f(x,th)}
f(x,\theta)=\log ( p(x|\theta_0)/ p(x|\theta)),
\end{align}
which is equivalent to
\[
p(x|\theta)=p(x|\theta_0)\exp(-f(x,\theta)). 
\]
We define $K(\theta)$ and $K_n(\theta)$ by
\begin{align}
K(\theta)&=\EE_X[f(X,\theta)|q],\label{eq:K(th)}
\\
K_n(\theta)&=\frac{1}{n}\sum_{i=1}^n f(X_i,\theta).\label{eq:Kn(th)}
\end{align}
Then $K(\theta)$ is a nonnegative function by the definition and 
\begin{align}
L(\theta)&=L(\theta_0)+K(\theta),
\label{eq:L1}
\\
L_n(\theta)&=L_n(\theta_0)+K_n(\theta),
\label{eq:Ln1}
\end{align}
resulting that 
\begin{align}
\frac{
\Omega(\theta) }{
\exp(-nL_n(\theta_0))
}
&=\exp(-nK_n(\theta)).
\end{align}
The set of all optimal parameters 
is equal to the set of all zero points of an analytic function
$K(\theta)$, 
\[
\Theta_0=\{\theta \;;\;K(\theta)=0\},
\]
which contains singularities in general. It has been difficult to study 
statistics and machine learning on the original parameter space because of singularities. 
The following theorem 
is the algebro-geometric foundation on which universal statistical theory can 
be constructed. 
\vskip3mm\noindent
{\bf Hironaka Resolution Theorem} \cite{Hironaka1964,Atiyah1970,Kashiwara1976}. 
There exist both 
a compact subset ${\cal M}$ of a $d$-dimensional analytic manifold   and 
a proper analytic function from ${\cal M}$ to $\Theta$ 
\[
g:{\cal M}\ni u \mapsto g(u)\in \Theta
\]
such that, in each local coordinate of ${\cal M}$, 
$K(g(u))$ is normal crossing, 
\begin{align}
K(g(u))&=u^{2k}\equiv u_1^{2k_1}u_2^{2k_2}\cdots u_d^{2k_d},
\label{eq:res01}
\\
\pi(g(u))|g'(u)|&=b(u)|u^h| \equiv b(u)|u_1^{h_1}u_2^{h_2}\cdots u_d^{h_d}|, 
\label{eq:res02}
\end{align}
where $k=(k_1,k_2,...,k_d)$ and $h=(h_1,h_2,...,h_d)$ are 
multi-indices of nonnegative integers, in which 
at least one $k_i$ is a positive integer. Here 
$b(u)>0$ is a positive analytic function, and $|g'(u)|$ is the
absolute value of the Jacobian determinant of $\theta=g(u)$. 
The correspondence between $\Theta\setminus\Theta_0$ and 
$g^{-1}(\Theta\setminus\Theta_0)$ is one-to-one in any 
neighborhood of $\Theta_0$. 
Note that a function $\theta=g(u)$ is called proper 
if the inverse image of a compact set is also compact.
\vskip3mm\noindent
This is the basic and most important theorem in algebraic geometry. 
There exists an algebraic algorithm by which both 
${\cal M}$ and $\theta=g(u)$ 
can be found by finite recursive blow-ups \cite{Hironaka1964} .
 If Newton diagram of $K(\theta)$
 is nondegerate, they are found by a toric modification, which 
was applied to statistics and machine learning \cite{Yamazaki2010}. 
Even if $K(\theta)$ is not an analytic function, if $K(\theta)=K_1(\theta)K_2(\theta)$ 
where $K_1(\theta)$ is analytic and $K_2(\theta)>0$, then the same theory can be 
derived in Bayesian statistics.  
For concrete examples of this theorem in statistics and machine learning, see \cite{Watanabe2009}. 
\vskip3mm\noindent
Based on Hironaka resolution theorem, 
Bayesian statistics of a pair on $\Theta$
\[
(p(x|\theta),\pi(\theta)) 
\]
is equivalent to that of a pair on ${\cal M}$ 
\[
(p(x|g(u)),\pi(g(u))|g'(u)|).
\]
It follows that
\begin{align}
\frac{\Omega(\theta) d\theta}{
\exp(-nL_n(\theta_0))}
&=\exp(-nK_n(g(u)))b(u)|u^h|du. 
\end{align}
Using the resolution theorem  eqs. (\ref{eq:res01}) and (\ref{eq:res02}), 
the {\it real log canonical threshold} (RLCT) $\lambda$ and 
{\it multiplicity} $m$ are defined by 
\begin{align}
\lambda&=\min_{L.C.}\min_{1\leq j\leq d}\left(\frac{h_j+1}{2k_j}\right),
\\
m&=\max_{L.C.}\#\left\{j;\frac{h_j+1}{2k_j}=\lambda\right\},
\end{align}
where $\min_{L.C.}$ and $\max_{L.C.}$
show the minimum and maximum values over all local coordinates respectively.
Here we define $(h_j+1)/(2k_j) =\infty$ for $k_j=0$,
and $\#$ means the number of elements of a set. 
Hence $0<\lambda<\infty$ and $1\leq m\leq d$. 
Since $\Theta$ is compact and $w=g(u)$ is proper, the number of all 
local coordinates is finite and the integration over $u$ is given by
the finite sum of integrations of local coordinates. Without loss of generality,
each local coordinate can be chosen as $[0,1]^d$ by the appropriate 
preparation of local parameter. 
For a given function $K(\theta)$, there are infinitely many pairs 
${\cal M}$ and $w=g(u)$ that give the resolution of singularities, 
however, neither RLCT nor its multiplicity depends on the
choice of such pairs, in other words, they are birational invariants
 \cite{Watanabe2009}. 

It is well known that 
the concept of the log canonical threshold plays an important
role in higher dimensional algebraic geometry \cite{Kollar1997} and 
the real log canonical threshold is the same concept in the
real algebraic geometry \cite{Saito2007}. 
It was found in \cite{Watanabe1999} that 
RLCT determines the accuracy of Bayesian statistics and machine learning. 
In many statistical models and learning machines, the common singularities 
often appear which is called Vandermonde type singularities \cite{Aoyagi2012}.
A concrete method and application to statistics and machine learning are 
introduced in \cite{Watanabe2009,Watanabe2018}. 

RLCT is determined uniquely for a given pair $(K(\theta),\pi(\theta))$. 
There are several mathematical properties. 
\begin{enumerate}
\item
If  $q(x)$ is regular for $p(x|\theta)$ and $\pi(\theta_0)>0$, then 
$\lambda=d/2$, $m=1$.
\item
If there exists $\theta_0\in\Theta_0$ such that $\nabla^2 L(\theta_0)=0$ and 
$\pi(\theta_0)>0$, then $0<\lambda < d/2$.
\item
Note that Jeffreys' prior is equal to zero at singularities, and if 
Jeffreys' prior is employed in singular models, then $\lambda\geq d/2$. 
\item
Assume that $\lambda_j$ $(j=1,2)$ are RLCTs of $(K_j(\theta_j),\pi_j(\theta_j))$.
Then 
\begin{itemize}
\item
RLCT of $(\sum_j K_j(\theta_j),\prod_j\pi_j(\theta_j))$ 
is equal to $(\sum_j\lambda_j)$. 
\item
RLCT of 
$(\prod_j K_j(\theta_j),\prod_j\pi_j(\theta_j))$ 
is equal to $(\min_j\lambda_j)$. 
\end{itemize}
\item
Assume that $\lambda_j$ $(j=1,2)$ are RLCTs of $(K_j(\theta),\pi_j(\theta))$
and that $K_1(\theta)\leq c_1K_2(\theta)$ 
and  $\pi_1(\theta)\geq c_2\pi_2(\theta)$ 
for some $c_1,c_2>0$. Then $\lambda_1\leq \lambda_2$. 

\end{enumerate}

These properties are helpful to study RCLTs of statistical models and learning
machines. In fact, 
RLCTs of important 
statistical models and learning machines were found by developing 
resolution procedures in 
neural networks \cite{Watanabe2001a,Aoyagi2012}, 
normal mixtures \cite{Yamazaki2003},
Poisson mixtures \cite{Sato2019}, multinomial mixtures \cite{WatanabeT2022}, 
general and nonnegative matrix facorizations \cite{Aoyagi2005,Hayashi2017},
Boltzmann machines \cite{Yamazaki2005a,Aoyagi2010}, 
hidden and general Markov models \cite{Yamazaki2005,Zwie2011}
and latent Dirichlet allocations \cite{Hayashi2021}.
Note that singularities in statistical models and learning machines make
the free energy and the generalization loss smaller 
if Bayesian inferences are employed, hence almost all learning machines
are singular \cite{Watanabe2007} and that's good \cite{Wei2022}. 
In mixture models, a Dirichlet distribution is often chosen for a prior 
distribution of the mixture ratio. Then the posterior distribution has a 
phase transition according to the hyperparameter of Dirichlet distribution
\cite{Watanabe2018,WatanabeT2022}. In the different phases, RCLTs are
different and the supports of asymptotic posterior distributions are different
\cite{Watanabe2018}.

\section{Two Mathematical Solutions}\label{section:two}

In this section, we show two mathematical problems in Bayesian
statistics are solved based on resolution theorem. 

\subsection{Singular Schwartz Distribution}\label{subsection:delta}

In this subsection, 
the first mathematical problem is solved on the 
resolution theorem. A method 
 how to analyze the set 
 $\{\theta\in\Theta\;;\;K(\theta)<\epsilon\}$ 
as $\epsilon\rightarrow +0$ is constructed even 
when an analytic set  $K(\theta)=0$  contains singularities. 
The mathematical method explained in this subsection is
based on the researches of singular Schwartz distribution by
 Gel'fand and Atiyah \cite{Atiyah1970,Kashiwara1976}. 

By the definition of RLCT, without loss of generality, we can assume that
\[
u=(u_a,u_b)\in\RR^m\times\RR^{d-m},
\]
which satisfies 
\begin{align}
\left(\frac{h_j+1}{2k_j}\right)&=\lambda\;\;\;(1\leq j\leq m), \nonumber
\\
\left(\frac{h_j+1}{2k_j}\right)&>\lambda\;\;\;(m+1\leq j\leq d), \nonumber
\end{align}
where $\lambda$ is RLCT. 
In other words, $u_a$ is the part of $u$ which gives $\lambda$ with
multiplicity $m$ and 
$u_b$ is the other part of $u$ which does not. 
A multi-index $\mu\in\RR^{d-m}$ is defined by
\[
\mu=\{\mu_j=-2\lambda k_j+h_j\;;\; m+1\leq j\leq d\}\in\RR^{d-m}.
\]
The zeta function $\zeta(z)$ for $z\in\CC$ and 
the state density function $v(t)$ for $t\in\RR$ are defined respectively by
\begin{align}
\zeta(z)&=
\int K(\theta)^z \pi(\theta)d\theta, 
\nonumber
\\
v(t)&=\int \delta(t-K(\theta))\pi(\theta)d\theta,
\end{align}
where $\zeta(z)$ in  $\Re(z) > -\lambda$ is well-defined 
by the integration over $\theta$, which 
can be analytically continued to the unique meromorphic function 
on the entire complex plane   \cite{Watanabe2009}. 
Then the zeta function is equal to the 
Mellin transform of the state density function, 
\[
\zeta(z)=\int v(t)\;t^{z}\;dt,
\]
therefore, the state density function is equal to
the inverse Mellin transform of the zeta function. 
By this correspondence, concrete calculation \cite{Watanabe2009}
 shows that 
\[
\frac{1}{(z+\lambda)^m}\Longleftrightarrow 
c_0 t^{\lambda-1}(-\log t)^{m-1},
\]
where $c_0$ is a constant. 
By using the resolution map $\theta=g(u)$, the zeta function is equal to the 
finite  sum  of the
integrations over $u\in[0,1]^d$ in local coordinates, 
\begin{align}
\zeta(z)& = \sum_{L.C.}\int_{[0,1]^d} K(g(u))^z \pi(g(u))|g'(u)|du
\nonumber
\\
&= \sum_{L.C.}\int_{[0,1]^d} u^{2kz}|u^h|b(u)du
\nonumber
\\
&= \sum_{L.C.}\int_{[0,1]^d} u^{2kz+h}b(0,u_b)du
\nonumber
\\
&+
\sum_{L.C.}\int_{[0,1]^d} u^{2kz+h}(b(u_a,u_b)-b(0,u_b))du.
\label{eq:2kz+h}
\end{align}
Let the first and second terms of eq.(\ref{eq:2kz+h}) 
be $\zeta_1(z)$ and $\zeta_2(z)$
respectively. 
The largest pole of $\zeta_1(z)$ 
is $(-\lambda)$ with the order $m$, and the 
largest pole of $\zeta_2(z)$ is smaller than $(-\lambda)$ or
its order is larger than $m$, because $b(u_a,u_b)-b(0,u_b)$ is divisible by 
$u_a$. By using the inverse Mellin transform, 
the following asymptotic expansion of the state density function
as $t\rightarrow +0$ is derived 
\cite{Watanabe2009,Watanabe2018},
\begin{align}
\delta(t-u^{2k})|u^h|b(u)du &=t^{\lambda-1}(-\log t)^{m-1}du^*
\nonumber
\\
&
+o(t^{\lambda-1}(-\log t)^{m-1}),
\label{eq:delta(t-u)}
\end{align}
where $du^*$ is an integration, 
\[
du^*=\frac{\delta(u_a)(u_b)^\mu\;b(u)}{2^m(m-1)! \prod_{j=1}^mk_j} \;du. 
\]
The asymptotic expansion eq.(\ref{eq:delta(t-u)}) will be employed in
the Bayesian theory in the following sections. 
Moreover, ${\rm Vol}(\varepsilon)$, which is the volume of 
the set of almost optimal parameters measured by the 
prior distribution, is given by
\[
{\rm Vol}(\varepsilon)
=\int_{K(\theta)<\varepsilon}
d\pi(\theta)
=\int_0^\varepsilon dt \int \delta(t-K(\theta))\pi(\theta)d\theta,
\]
is given by 
\begin{align}
{\rm Vol}(\varepsilon)&=\sum_{L.C.}\int_0^{\epsilon}dt\int \delta(t-u^{2k})|u^h|b(u)du
\nonumber
\\
&\propto \epsilon^{\lambda}(-\log \epsilon)^{m-1}+\mbox{small order}.
\end{align}
It follows that 
\[
\lambda=\lim_{\varepsilon\rightarrow +0}
\frac{\log{\rm Vol}(\varepsilon)}{\log\varepsilon}.
\]
which shows that RLCT can be understood as the generalized
 dimension of the set  $K(\theta)=0$.  As is shown in the following 
sections, RLCT determines the accuracy of Bayesian inference.

\subsection{Empirical Process and Renormalized Posterior}\label{subsection:empirical}

In this subsection, 
the second mathematical problem is solved on the 
resolution theorem. 
It has been difficult to treat the empirical process
on the original parameter space 
in the neighborhood of singularities of $K(\theta)=0$. 

The function $K_n(\theta)$ in eq.(\ref{eq:Kn(th)}) is rewritten as 
\begin{align}
K_n(\theta)&=K(\theta)-\frac{1}{\sqrt{n}}K(\theta)^{1/2}
\xi_n(\theta),
\end{align}
where 
\[
\xi_n(\theta)=\frac{1}{\sqrt{n}}\sum_{i=1}^n
\left\{
\frac{
K(\theta) - f(X_i,\theta)
}{\sqrt{K(\theta)}}
\right\}
\] 
which is not a well-defined function in general at $K(\theta)=0$. 

A function $f(x,\theta)$ is said to have a relatively finite 
variance, if  there exists a constant $c_0>0$ such that,
for an arbitrary $\theta\in\Theta$,
\begin{align}\label{eq:relbd}
\int q(x)f(x,\theta)^2dx\leq c_0 K(\theta). 
\end{align}
Since $\Theta$ is a compact set, this inequality holds in 
$K(\theta)>0$. If $q(x)$ is regular for $p(x|\theta)$, then 
$K(\theta)$ is a positive definite quadratic among $K(\theta)=0$, 
hence eq.(\ref{eq:relbd}) holds. If $q(x)$ is realizable by $p(x|\theta)$, then
also eq.(\ref{eq:relbd}) holds. In this paper, we study the case when 
eq.(\ref{eq:relbd}) holds. When $q(x)$ is unrealizable by and singular
for $p(x|\theta)$, then eq.(\ref{eq:relbd}) does not hod in general, 
resulting that the free energy has the different behavior from this paper \cite{Watanabe2010a,Nagayasu2022}. 

Assume that 
$f(x,\theta)$ has a relatively finite variance and that $K(g(u))=u^{2k}$ is normal 
crossing. Then by the factor theorem, for each $u_j$, 
 $f(x,g(u))^2$ is divisible by $u_j^{2k}$, hence $f(x,g(u))^2$ is divisible by 
$K(g(u))$. That is say, 
there exists a function $a(x,u)$ which is analytic for $u$ and
\[
f(x,g(u))=a(x,u)u^k.
\]
Then by the definition,  $\EE_X[a(X,u)|q]=u^k$. It follows that 
$\xi_n(u)\equiv\xi_n(g(u))$ is given by 
\begin{align}
\xi_n(u)&=
\frac{1}{\sqrt{n}}\sum_{i=1}^n
\{\EE_X[a(X,u)|q]-a(X_i,u)\}
\\
&=\frac{1}{\sqrt{n}}\sum_{i=1}^n \{   u^k - a(X_i,u)\} 
\end{align}
 is a well-defined function of $u$ and
\begin{align}
K_n(g(u))&=u^{2k}-
\frac{1}{\sqrt{n}}u^k\xi_n(u). 
\end{align}
Here $\xi_n(u)$ 
is an empirical process which converges to a Gaussian process 
in distribution $\xi_n(u)\rightarrow \xi(u)$ on each local coordinate $[0,1]^d$. 
The Gaussian process $\xi(u)$ 
is uniquely characterized by the conditions, 
\begin{align}
\EE[\xi_n(u)|q]&=\EE_\xi[\xi(u)]=0,
\\
\EE[\xi_n(u)\xi_n(v)|q]&=\EE_\xi[\xi(u)\xi(v)]=\EE_X[a(X,u)a(X,v)|q]-u^kv^k,
\end{align}
where $\EE_\xi[\;\;]$ shows the expectation value over $\xi(u)$. Note that 
and $\xi(u)$ and $\xi_n(u)$ have
the same expectation and covariance. 
Since the convergence of $\xi_n(u)\rightarrow \xi(u)$ is 
given by the uniform topology on $u\in{\cal M}$, it follows \cite{Watanabe2009} that 
\[
\lim_{n\rightarrow\infty}
\EE[\sup_{u}|\xi_n(u)|^2|q]=\EE_\xi[\sup_{u}|\xi(u)|^2]\]
where $\EE_\xi[\;\;]$ means the average over Gaussian process $\xi(u)$.
\vskip3mm\noindent
{\bf Standard Representation}. 
Let us summarize the two mathematical solutions
 in subsections \ref{subsection:delta} and \ref{subsection:empirical}.
Even for nonidentifable and singular statistical models and learning 
machines, there exists an appropriate representation of a parameter
$\theta=g(u)$, such that 
the average and
the log likelihood functions eq.(\ref{eq:Ln}) and eq.(\ref{eq:L}) 
can be rewritten as 
\begin{align}
L(g(u))&=L(\theta_0)+u^{2k}
\label{eq:standardL}
\\
L_n(g(u))&=L_n(\theta_0)+u^{2k}-\frac{1}{\sqrt{n}}u^k\xi_n(u),
\label{eq:standard}
\end{align}
in each local coordinate of the manifold ${\cal M}=g^{-1}(\Theta)$. 
In other words, for an arbitrary triple $(q(x),p(x|\theta),\pi(\theta))$,
universal properties of Bayesian statistics can be derived 
from the normal crossing representation eq.(\ref{eq:standard}). 
\vskip3mm\noindent
In this paper, we mainly study Bayesian statistics, however, 
the standard form eq.(\ref{eq:standard}) can be applied to 
the maximum likelihood and maximum {\it a posterior} methods 
\cite{Watanabe2009}. For example, let $\hat{\theta}$ be 
the maximum likelihood estimator. If there exists an odd number 
in $k=(k_1,k_2,...,k_d)$, then by the compactness of $\Theta$, 
\begin{align}
L(\hat{\theta})&=L(\theta_0)+\frac{1}{2n}\sup_{g(u)\in\Theta_0}|\xi_n(u)|^2,
\\
L_n(\hat{\theta})&=L_n(\theta_0)-\frac{1}{2n}\sup_{g(u)\in\Theta_0}|\xi_n(u)|^2,
\end{align}
where $\displaystyle \sup_{g(u)\in\Theta_0}$ shows the
supremum value over all optimal parameters $\Theta_0$. 
The function $L_n(\hat{\theta})$ 
is made smaller by the maximum likelihood method 
than by Bayesian method, however,
the function $L(\hat{\theta})$ becomes larger. 
This is a reason why Bayesian inference is better than the maximum
likelihood method in nonidentifiable or singular models.

\section{Three Applications to Statistics} \label{section:three}

In this section, we show that 
two mathematical solutions
 in subsections \ref{subsection:delta} and \ref{subsection:empirical}
are useful in three applications to statistics. First, 
the the posterior distribution is represented by a well-defined 
renormalized posterior distribution in subsection
\ref{subsection:111},  second,
the asymptotic behavior of the free energy and its estimating methods are 
derived in subsection \ref{subsection:222},  lastly, 
the universal formula among the generalization loss, the cross validation loss,
and the information criterion are established in subsection \ref{subsection:333}.

\subsection{Renormalized Posterior Distribution}\label{subsection:111}

First, the asymptotic behavior of the posterior distribution is represented by
using the renormalized posterior distribution. 

By applying eqs. (\ref{eq:delta(t-u)}) and (\ref{eq:standard}), 
we can derive the asymptotic behavior of singular posterior distribution. 
\begin{align}
\frac{\Omega(\theta) d\theta}{\exp(-nL_n(\theta_0))}
&=\exp(-nK_n(\theta))\pi(\theta)d\theta \nonumber 
\\
&=\exp(-nu^{2k}+\sqrt{n}\;u^k \xi_n(u))b(u)|u^h|\;du \nonumber
\\
&=\int dt\; \delta(t-nu^{2k})\exp(-t+\sqrt{t} \xi_n(u))b(u)|u^h|\;du
\nonumber
\\
&= \frac{(\log n)^{m-1}}{n^{\lambda}}
\int dt\; t^{\lambda-1}
\exp(-t+\sqrt{t} \xi_n(u))\;du^*
\nonumber
\\
&
+ o_p\left(\frac{(\log n)^{m-1}}{n^{\lambda}}\right). 
\label{eq:omega}
\end{align}
The renormalized posterior distribution is defined by its expectation value
$\langle F(t,u) \rangle$ 
of a given 
function $F(t,u)$, 
\begin{align}\label{eq:F(t,u)}
\langle F(t,u)\rangle
=\frac{\sum_{L.C.} \int du^* \;dt\; F(t,u)\; t^{\lambda-1}\;\exp(-t+\sqrt{t}\;\xi_n(u))
}{
\sum_{L.C.} \int du^* \;dt\; t^{\lambda-1}\;\exp(-t+\sqrt{t}\;\xi_n(u))
}. 
\end{align}
Also $\langle F(t,u) \rangle_{\infty}$ is defined by eq.(\ref{eq:F(t,u)}) 
whose $\xi_n(u)$ is replaced by $\xi(u)$. 
Then the renormalized posterior distribution satisfies 
\[
\lim_{n\rightarrow\infty }
\EE[\langle F(t,u)\rangle|q] =
\EE_{\xi}[\langle F(t,u)\rangle_{\infty}].
\]
By using the renormalized posterior distribution, the asymptotic 
behavior of the posterior distribution is derived as follows. 
By the relation $f(x,g(u))=a(x,u)u^k$ and  $t=n\;u^{2k}$,  
the correspondence between the posterior distribution and 
the renormalized posterior distribution is derived for $\alpha\geq 0 $, 
\begin{align}\label{eq:renormalize1}
\EE_\theta[ f(x,\theta)^\alpha ]
=\frac{1}{n^{\alpha/2}}\langle (a(x,u)\sqrt{t})^\alpha\rangle
+o_p(1/n^{\alpha/2}). 
\end{align}
This equation  shows that the log density ration function bas the
$(1/n^{1/2})$ order. 
Moreover, there exist two equations which hold by the renormalized
posterior distribution. First,
by using  a partial integration over $t$, it follows that 
\begin{align}\label{eq:renormalize2}
\langle t^{\alpha}\rangle&=(\lambda+\alpha-1)\langle t^{\alpha-1}\rangle
+\frac{1}{2}
\langle t^{\alpha-1/2}\;\xi_n(u)\rangle .
\end{align}
 Second,
let us define two random variables $\langle \sqrt{t}\xi_n(u)\rangle$ and 
\[
V(\xi_n)=\EE_X[ \langle t a(X,u)^2\rangle -\langle \sqrt{t}a(X,u)\rangle^2|q].
\]
Then by using the partial integration over the functional space $\xi(u)$, it follows that
\[
\EE_\xi[\langle \sqrt{t}\;\xi(u)\rangle_{\infty}]=\EE_\xi[V(\xi)].
\]
By using these properties,  the asymptotic values of 
expectations by the posterior distribution are represented by the 
renormalized posterior distributions. 

\subsection{Asymptotic Free Energy}\label{subsection:222}

Second, the asymptotic free energy, which is equal to the 
minus log marginal likelihood, is derived. 
By applying eq.(\ref{eq:omega}) to eq.(\ref{eq:FOmega}), it follows that 
\begin{align}
F_n&=nL_n(\theta_0)+\lambda\log n-(m-1)\log\log n 
\nonumber
\\
& +\chi(\xi_n) +o_p(1),
\label{eq:Fnrlct1}
\\
\EE[F_n|q]&=nL(\theta_0)+\lambda\log n-(m-1)\log\log n 
\nonumber
\\
&+\EE_\xi[\chi(\xi)]+o(1),
\label{eq:Fnrlct2}
\end{align}
where $\chi(\xi_n)$ is defined by 
\[
\chi(\xi_n)=-\log \left(
\sum_{L.C.}\int dt\; t^{\lambda-1}\exp(-t+\sqrt{t}\;\xi_n(u))du^*
\right),
\]
which converges to $\chi(\xi)$  in distribution. 

In general, 
it needs heavy computational costs
 to calculate the free energy $F_n$, hence alternative methods are desired. 
The real log canonical threshold depends on not only a pair  $p(x|\theta)$ 
and $\pi(\theta)$ but also an unknown probability distribution $q(x)$. 
Hence neither eqs. (\ref{eq:Fnrlct1}) nor (\ref{eq:Fnrlct2}) can be directly 
applied  to numerical calculation of the free energy. 
Two methods were proposed for solving this problem for the case when
the posterior distribution cannot be approximated by any normal distribution. 
The first method was proposed by \cite{Drton2017} that,
by using the estimated RLCT $\hat{\lambda}$, the singular BIC 
\[
{\rm sBIC}=nL_n(\hat{\theta})+\hat{\lambda}\log n
\]
is defined, where $\hat{\theta}$ is 
the maximum likelihood estimator. This is a generalized version of 
BIC \cite{Schwarz1978} of regular statistical models to general
singular models. 

The second method \cite{Watanabe2013} is as follows. 
We introduce a partition function for an inverse temperature 
$\beta>0$, 
\[
{\cal F}(\beta)=-\log \int \prod_{i=1}^np(X_i|\theta)^{\beta}
\pi(\theta)d\theta
\]
Then ${\cal F}(0)=0$ and ${\cal F}(1)=F_n$, resulting that
there exists $0<\beta^*<1$ such that
\begin{align}
F_n&={\cal F}'(\beta^*)
\nonumber
\\
&= \frac{\int nL_n(\theta)
\prod_{i=1}^np(X_i|\theta)^{\beta^*}
\pi(\theta)d\theta
}{\int 
\prod_{i=1}^np(X_i|\theta)^{\beta^*}
\pi(\theta)d\theta
}
\label{eq:WBIC}
\\
&= nL_n(\theta_0)+
\frac{
\int nK_n(\theta)
\prod_{i=1}^n\exp(-n\beta^*K_n(\theta))
\pi(\theta)d\theta
}{
\int 
\prod_{i=1}^n\exp(-n\beta^*K_n(\theta))
\pi(\theta)d\theta
}
\label{eq:beta*}.
\end{align}
Let the second term of the right hand side of eq.(\ref{eq:beta*}) be
$R_n$. Then  by the same way as the renormalized
posterior distribution is derived, 
\[
R_n
=\frac{\sum_{L.C.} \int du^* \;dt\;
(t-\sqrt{t}\xi_n(u))
\; t^{\lambda-1}\;\exp(-\beta^*t+\beta^*\sqrt{t}\;\xi_n(u))
}{
\sum_{L.C.} \int du^* \;dt\; t^{\lambda-1}\;\exp(-\beta^*t+\beta^*\sqrt{t}\;\xi_n(u))
}. 
\]
In this integration, by setting $\beta^*=1/\log n$ and by replacing 
$t$ and $dt$ by $t/\beta^*$ and 
$dt/\beta^*$ in the integration respectively, it follows that 
\[
R_n=\lambda \log n+o_p(\log n).
\]
Therefore by the definition 
\[
{\rm WBIC }={\cal F}'(1/\log n),
\] 
we obtain 
\[
{\rm WBIC}
=nL_n(\theta_0)+\lambda\log n
+o_p(\log n). 
\]
has the same asymptotic expansion as $F_n$ according to
the order $\log n$. In the numerical calculation of WBIC, the posterior 
distribution with the inverse temperature $1/\log n$ is necessary. 
An efficient algorithm to generate such posterior distribution in mixtures 
models are proposed \cite{Watanabe2021}. 

In sBIC, no averaging calculation on the parameter set is required, 
but the theoretical results about RLCTs for several models are necessary. 
In
WBIC, averaging on the parameter set is necessary, but the theoretical
results about RLCTs are not required.

In the variational Bayes approaches, the Bayesian posterior distribution is
approximated by an independent distribution $r(\theta)=r_(\theta_1)r_2(\theta_2)$, 
and the variational free energy is defined as the minimization of the functional 
\[
F_{vb}=\inf_{r=r_1\cdot r_2}\left\{
-S(r)- \int r(\theta)\log \Omega(\theta)d\theta
\right\}
\]
where $S(r)$ is the entropy of $r(\theta)$ and 
$\displaystyle \inf_{r=r_1\cdot r_2}$ is the infimum value overall probability distributions
that are represented by  
$r(\theta)=r_1(\theta_1)r_2(\theta_2)$. In singular cases,
the variational free energy has a different coefficient of the 
$\log n$ term from that of Bayesian free energy 
 \cite{WatanabeK2006,Yamazaki2013,Nakajima2019,Kariya2020a}.
The Kullback-Leibler divergence between the posterior and 
approximated distributions is equal to 
the difference between the Bayes and variational free energies, 
hence the accuracy of the variational approximation is clarified by examining 
both free energies. 

\subsection{Generalization Loss and Its Estimators} 
\label{subsection:333}

Third, we show that the asymptotic behaviors of generalization loss 
and its estimators are clarified. 

Let us introduce the functional cumulant generating functions \cite{Watanabe2010}
for $\alpha\in\RR$, 
\begin{align}
{\cal G}(\alpha)&=\EE_X[\log\EE_\theta[p(X|\theta)^{\alpha}]|q],
\\
{\cal T}(\alpha)&=\frac{1}{n}\sum_{i=1}^n\log
\EE_\theta[p(X_i|\theta)^{\alpha}],
\end{align}
by which 
random variables $G_n$, $T_n$, $C_n$ and $W_n$ are
represented, 
\begin{align}
G_n&=-{\cal G}(1),
\\
T_n&=-{\cal T}(1),
\\
C_n&={\cal T}(-1),
\\
W_n&=-{\cal T}(1)+{\cal T}''(0).
\end{align}
Note that ${\cal G}(0)={\cal T}(0)=0$. 
It is shown in \cite{Watanabe2010} that
$\EE[{\cal G}(\alpha)|q]=\EE[{\cal T}(\alpha)|q]$ and that, 
for $k\geq 2$, 
\begin{align}
{\cal G}^{(k)}(\alpha)&=O_p(1/n^{k/2}),
\\
{\cal T}^{(k)}(\alpha)&=O_p(1/n^{k/2}).
\end{align}
The four random variables are represented by the functional cumulant 
generating functions, 
\begin{align}
G_n&=-{\cal G}'(0)-\frac{1}{2}{\cal G}''(0)+O_p(1/n^{3/2}),
\\
T_n&=-{\cal T}'(0)-\frac{1}{2}{\cal T}''(0)+O_p(1/n^{3/2}),
\\
C_n&=-{\cal T}'(0)+\frac{1}{2}{\cal T}''(0)+O_p(1/n^{3/2}),
\\
W_n&=-{\cal T}'(0)+\frac{1}{2}{\cal T}''(0)+O_p(1/n^{3/2}),
\end{align}
By the definition of the log density ratio function eq.(\ref{eq:f(x,th)}), 
it follows that 
\begin{align}
-{\cal G}'(0)&=L(\theta_0)+\EE_X[\EE_\theta[f(X|\theta)]|q],
\\
{\cal G}''(0)&=\EE_X[\VV_\theta[f(X|\theta)]|q],
\\
-{\cal T}'(0)&=L_n(\theta_0)+
\frac{1}{n}\sum_{i=1}^n \EE_\theta[f(X_i|\theta)],
\\
{\cal T}''(0)&=\frac{1}{n}\sum_{i=1}^n \VV_\theta[f(X_i|\theta)].
\end{align}
By applying eq.(\ref{eq:renormalize1}) and eq.(\ref{eq:renormalize2}),
\begin{align}
G_n&=L(\theta_0)+\frac{1}{n}
\left(\lambda+\frac{1}{2}\langle \sqrt{t}\;\xi_n(u)\rangle
-\frac{1}{2}V(\xi_n)
\right)+o_p\left(\frac{1}{n}\right),\label{eq:Gn3}
\\
T_n&=L_n(\theta_0)+\frac{1}{n}
\left(\lambda-\frac{1}{2}\langle \sqrt{t}\;\xi_n(u)\rangle
-\frac{1}{2}V(\xi_n)
\right)+o_p\left(\frac{1}{n}\right),\label{eq:Tn3}
\\
C_n&=L_n(\theta_0)+\frac{1}{n}
\left(\lambda-\frac{1}{2}\langle \sqrt{t}\;\xi_n(u)\rangle
+\frac{1}{2}V(\xi_n)
\right)+o_p\left(\frac{1}{n}\right),\label{eq:Cn3}
\\
W_n&=L_n(\theta_0)+\frac{1}{n}
\left(\lambda-\frac{1}{2}\langle \sqrt{t}\;\xi_n(u)\rangle
+\frac{1}{2}V(\xi_n)
\right)+o_p\left(\frac{1}{n}\right).\label{eq:Wn3}
\end{align}
Let $\nu=\EE_\xi[V(\xi)]/2$ be the {\it singular fluctuation}. It follows that 
\begin{align}
\EE[G_n|q]&=L(\theta_0)+\frac{\lambda}{n}+
o\left(\frac{1}{n}\right),
\\
\EE[T_n|q]&=L(\theta_0)+\frac{\lambda-2\nu}{n}+o\left(\frac{1}{n}\right),
\\
\EE[C_n|q]&=L(\theta_0)+\frac{\lambda}{n}+
o\left(\frac{1}{n}\right),
\\
\EE[W_n|q]&=L(\theta_0)+\frac{\lambda}{n}+
o\left(\frac{1}{n}\right).
\end{align}
These results theoretically clarified the asymptotic behaviors of 
the generalization loss and its estimators, on which we can make 
Bayesian model evaluation methods. 
The leave-one-out cross validation (LOOCV) and WAIC can be employed even if
$q(x)$ is singular for $p(x|w)$, whereas neither AIC \cite{Akaike1974} nor DIC \cite{Spiegel2002}. 
When a leverage sample point is contained in a sample, 
the importance sampling cross validation eq.(\ref{eq:Cn2}) becomes unstable
\cite{Peruggia1997,Epifani2008}, and 
the difference between LOOCV and WAIC 
becomes larger \cite{Watanabe2018,Watanabe2021a}. 
The improved version of numerical calculation of the cross validation was proposed in \cite{Vehtari2017}. 

If a sample is independent, LOOCV and 
WAIC are equivalent to each other. However, if otherwise, they may not
be equivalent. For example, 
in regression problems where the conditional probability distribution 
$q(y|x)$ of an output $Y$ for a given input $X$, 
the input samples $\{X_i\}$ may dependent or fixed and $\{Y_i\}$ are
conditionally independent. 
In such cases,  LOOCV does not estimate the conditional generalization loss, 
 whereas WAIC does \cite{Watanabe2018,Watanabe2021a}.

From eqs. (\ref{eq:Gn3}), (\ref{eq:Tn3}), (\ref{eq:Cn3}), and (\ref{eq:Wn3}), 
it is derived that 
the leave-one-out cross validation
and the information criterion have the inverse correlation to the generalization loss, 
\cite{Watanabe2010}, 
\begin{align}
(G_n-L(\theta_0))+(C_n-L_n(\theta_0)) & =\frac{2\lambda}{n}+o_p(\frac{1}{n}),
\\
(G_n-L(\theta_0))+(W_n-L_n(\theta_0)) & =\frac{2\lambda}{n}+o_p(\frac{1}{n}). 
\end{align}
Although the cross validation and the information criterion are useful in many statistical
applications, these properties clarified a disadvantages of them. Improved methods, 
adjusted cross validation and information criteria, 
combining leave-one-out and hold-out cross validations have been proposed 
\cite{Watanabe2022a}, which make the variance of the estimators smaller. 

From the viewpoint of the bias and variance problem, the effect of singularities
are also studied. 
The probabilistic behaviors of the generalization losses when the optimal parameter 
$\Theta_0$ is in a neighborhood of singularities were also clarified \cite{Watanabe2001}.
In singular models and machines, there are phase transitions as sample size increases \cite{Watanabe2003}.

\section{Conclusion}\label{section:conc}

We have reviewed recent advances in the research filed of algebraic geometry
and Bayesian statistics. The two mathematical problems caused by singular 
log likelihood function were resolved by an algebro-geomegtric transform.  
There are two birational
invariants which determine Bayesian statistics. 
The former is the real log canonical
threshold which clarifies the singular dimension of a statistical model and
a prior distribution. The latter 
is the singular fluctuation which indicates the functional variance of the 
log likelihood function. Based on the theoretical properties of
 these two concepts, three statistical 
problems were overcome. First, the posterior distribution was represented by
 a renormalized posterior distribution 
defined on a manifold. Second, the asymptotic behavior of the free energy was
clarified and its estimation methods were constructed.  Lastly,  
universal formulas 
between generalization loss, cross validation, and information criterion were 
derived. These mathematical and statistical results are  now  being 
used in data science and
artificial intelligence.

\section*{ Data Availability and Conflict of Interest}

\subsection*{\it Data Availability}
 
Data sharing is not applicable to this article as no data sets were generated or analyzed during the current study.
 
\subsection*{\it Conflict of interest}

The corresponding author states that there is no conflict of interest.


\begin{thebibliography}{99}

\bibitem{Akaike1974}
Akaike, H. A new look at the statistical model identification, 
 IEEE Transactions on Automatic Control. Vol.19, No.6, pp.716-723, 1974.

\bibitem{Akaike1980}
Akaike, H.
On the transition of the paradigm of statistical inference.
The proceedings of the Institute of Statistical Mathematics,
Vol.27,	 pp.5-12. 1980. 

\bibitem{Amari2001}
Amari, S. 
Differential and algebraic geometry in multilayer perceptrons.
IEICE Transactions on Fundamentals, Vol.E84-A, pp.31-38, 2001.

\bibitem{Amari1992}
Amari, S., Fujita, N., Shinomoto,S. Four types of Leaning Curves.
Neural Computation, Vol.4, pp.605-618, 1992. 

\bibitem{Amari1993}
Amari, S., Murata, N. 
Statistical theory of learning curves under entropic loss criterion. Neural Computation, Vol.5, pp.140-153,
Nueral Computation. 

\bibitem{Aoyagi2005}
Aoyagi, M., Watanabe, S. Stochastic complexities of reduced rank regression in Bayesian estimation. Neural Networks, Vol.18, pp.924-933, 2005.

\bibitem{Aoyagi2010}
Aoyagi, M. Stochastic complexity and generalization error of a restricted
Boltzmann machine in Bayesian estimation
Journal of Machine Learning Research, Vol.11, pp.1243-1272, 2010. 

\bibitem{Aoyagi2012}
Aoyagi, M., Nagata, K. Learning coefficient of generalization error in Bayesian estimation and Vandermonde matrix type singularity. 
Neural Computation, vol. 24, No. 6, pp.1569-1610, 2012.

\bibitem{Atiyah1970}
 Atiyah. M. F. Resolution of singularities and division of distributions
communications on pure and applied mathematics. Vol.23, no.2. pp.145-150. 1970.

\bibitem{Binmore2017}
Binmore, K. On the foundations of decision theory. Homo Oeconomicus, Vol.34, pp.259
-273, 2017.

\bibitem{Box1976}
Box, G. E. P. Science and statistics. Journal of the American Statistics Association. Vol.71, pp.791-799, 1976. 

\bibitem{Drton2017}
Drton, M., Plummer, M. A Bayesian information criterion for singular models. J. R. Statist. Soc.  B., Vol.56, pp.1-38, 2017.

\bibitem{Epifani2008}
Epifani, I., MacEchern, S. N., Peruggia, M. Case-Deletion 
importance sampling estimators: Central limit 
theorems and related results. Electric Journal of 
Statistics, Vol.2, pp.774-806, 2008.

\bibitem{Fukumizu1996}
Fukumizu, K., A regularity condition of the information matrix of a multilayer perceptron network.
Neural Networks, Vol.9, pp.871-879, 1996.

\bibitem{Gelfand1992}
Gelfand, A. E., Dey, D. K., Chang, H. 
Model determination using predictive distributions with
  implementation via sampling-based method. 
Technical Report, Department of statistics,
Stanford University, Vol.462, pp. 147-167, 1992.

\bibitem{GelmanBDA} Gelman, A.,
Carlin, J. B.,
Stern H.S.,
Dunson, D.B.,
Vehtari, A.,
Rubin D.B.
Bayesian data analysis III. CRC Press. 2013. 

\bibitem{Gelman2013}
Gelman, A., Shalizi, C. S. 
Philosophy and the practice of Bayesian statistics. 
British Journal of Mathematical and Statistical Psychology. 66, pp.8-38, 2013.

\bibitem{Gelman2014}
Gelman, A., Hwang, J., Vehtari, A.
 Understanding predictive information criteria for Bayesian models.
Statistics and Computing, Vol.24, pp.997--1016, 2014.

\bibitem{Hagiwara1993}
Hagiwara, K., Toda, N., Usui, S. 
On the problem of applying AIC to determine the structure of a layered feedforward neural network.
Proc. of 1993 International Conference on Neural Networks, 
Vol. 3, pp.2263-2266, 1993.

\bibitem{Hartigan1985} 
Hartigan, J. A. A failure of likelihood asymptotics for normal mixtures. Proc. of
Berkeley Conference in Honor of J.Neyman and J.Kiefer, Vol.2, pp.807?810, 1985.


\bibitem{Hayashi2017}
Hayashi, N., Watanabe, S. 
Upper bound of Bayesian generalization error in non-negative matrix
 factorization, Neurocomputing 266, pp.21-28, 2017. 

\bibitem{Hayashi2021}
Hayashi, N. 
The exact asymptotic form of Bayesian generalization error in latent Dirichlet allocation. 
Neural Networks, Vol.137, pp.127-137, 2021.

\bibitem{Hironaka1964}
Hironaka, H.
Resolution of singularities of an algebraic variety over a field of characteristic zero. I,II.  Ann. of Math., Vol.79, pp.109-326, 1964. 

\bibitem{Kariya2020}
Kariya, N. Watanabe, S. 
Asymptotic analysis of singular likelihood ratio of normal mixture by Bayesian learning theory for testing homogeneity. 
Communications in Statistics-Theory and Methods, pp.1-18, Vol.51, 2020. 

\bibitem{Kariya2020a}
Kariya, N. Watanabe, S. 
Testing homogeneity for normal mixture models: variational Bayes approach. 
IEICE TRANSACTIONS on Fundamentals of Electronics, Communications and Computer Sciences. Vol.103, pp.1274-1282, 2020. 

\bibitem{Kashiwara1976}
Kashiwara, M. B-functions and holonomic systems. 
Rationality of roots of B-functions. 
Inventiones mathematicae, Vol.38, pp.33-53, 1976. 

\bibitem{Kollar1997}
 Koll\'{a}r, J.
Singularities of pairs, 
Proceedings of Symp. Pure Math., A.M.S. Vol.62, Part 1, pp. 221-287.
1997. 

\bibitem{McElreath2020}
McElreath, S.
Statistical Rethinking: A Bayesian course with examples in R and STAN. 2nd edition. CRC Press, 2020. 

\bibitem{Murata1995}
Murata, N., Yoshizawa, S., Amari, S. 
Network information criterion-determining the number of hidden units for an artificial neural network model.
IEEE transactions on neural networks, Vol.5 ,pp.865-872, 1995. 

\bibitem{Nagata2008}
Nagata, K. Watanabe, S. 
Asymptotic behavior of exchange ratio in exchange Monte Carlo method. 
 Neural Networks, Vol. 21, No. 7, pp. 980-988, 2008.
 
\bibitem{Nagayasu2022}
Nagayasu, S.,  Watanabe, S. 
 Asymptotic behavior of free energy when optimal probability distribution is not unique.
Neurocomputing,
Vol.500, pp.528-536, 2022. 

\bibitem{Nakajima2019}
Nakajima, S. Watanake, K. Sugiyama, M. Variational Bayesian Learning Theory, Cambridge University Press, 2019. 

\bibitem{Peruggia1997}
 Peruggia, M.
 On the variability of case-detection importance sampling weights in
  the Bayesian linear model.
Journal of American Statistical Association, Vol.92, 
pp.199-207, 1997.


\bibitem{Saito2007}
Saito, M. On real log canonical thresholds, 
arxiv:0707.2308, 2007. 

\bibitem{Sato2019}
Sato, K., Watanabe, S. 
Bayesian generalization error of Poisson mixture and simplex Vandermonde 
matrix type singularity. arXiv:1912.13289, 2019. 

\bibitem{Schwarz1978}
Schwarz, G. 
Estimating the dimension of a model. Vol. 6, Np.2 Annals of Statistics, pp.461-464. 1978. 

\bibitem{Spiegel2002}
Spiegelhalter, D. J., Best, N. G.,
Carlin, B. P., Linde, A. Bayesian measures of model complexity and
fit. Journal of Royal Statistical Society, Series B, Vol.64, No.4, pp.583-639, 2002. 


\bibitem{Vehtari2002}
Vehtari, A., Lampinen, J.
Bayesian model assessment and comparison using cross-validation
  predictive densities.
Neural Computation, Vol.14, no.10, pp.2439--2468, 2002.

\bibitem{Vehtari2017}
Vehtari, A., Gelman, A., Gabry, J. 
Practical Bayesian model evaluation using leave-one-out cross-validation and WAIC.
Statistics and computing. Vol. 27, No.5, pp.1413-1432, 2017. 

\bibitem{WatanabeK2006}
Watanabe, K., Watanabe, S.
Stochastic complexities of Gaussian mixtures in variational Bayesian approximation.
Journal of Machine Learning Research
Vol.7,pp.625-644, 2006. 

\bibitem{Watanabe1995}
Watanabe, S. A generalized Bayesian framework for neural networks with singular Fisher information matrices.
 Proc. of International Symposium on Nonlinear Theory and Its Applications, pp.207-210, 1995.

\bibitem{Watanabe1999}
S.Watanabe, S. Algebraic analysis for singular statistical estimation.
Lecture Notes in Computer Sciences, Vol.1720, pp.39-50, 1999.


\bibitem{Watanabe2001}
Watanabe, S. Algebraic geometrical methods for hierarchical learning machines.
Neural Networks, Vol.14, pp.1049-1060, 2001. 

\bibitem{Watanabe2001a}
Watanabe, S. 
Learning efficiency of redundant neural networks in Bayesian estimation. 
IEEE Transactions on Neural Networks.Vol.12, pp.1475-1486, 2001.

\bibitem{Watanabe2001b}
Watanabe, S. 
Algebraic analysis for nonidentifiable learning machines. 
Neural Computation. Vol.13, pp.899-933, 2001. 

\bibitem{Watanabe2003}
Watanabe, S., Amari, S. 
Learning coefficients of layered models when the true
distribution mismatches the singularities. 
Neural Computation, Vol.15, pp. 1013-1033, 2003. 

\bibitem{Watanabe2007}
Watanabe, S. 
Almost all learning machines are singular.
IEEE Symposium on Foundations of Computational Intelligence, 
pp.383-388, 2017. 

\bibitem{Watanabe2009} 
Watanabe, S. Algebraic geometry and statistical learning theory.
Cambridge University Press. 2009. 

\bibitem{Watanabe2010}
Watanabe, S. Asymptotic equivalence of Bayes cross validation and widely
applicable information criterion in singular learning theory. Journal of Machine Learning Research. Vol.11, pp.3571-3594, 2010.

\bibitem{Watanabe2010a}
Watanabe, S. 
Asymptotic learning curve and renormalizable condition in statistical learning theory, Journal of Physics Conference Series, Vol. 233, No. 1, 2010.

\bibitem{Watanabe2013}
Watanabe, S. 
A widely applicable Bayesian information criterion. Journal of Machine Learning Research. Vol. 14, pp.867-897, 2013. 

\bibitem{Watanabe2018} 
Watanabe, S. Mathematical theory of Bayesian statistics. 
CRC Press, 2018. 


\bibitem{Watanabe2018a}
Watanabe, S. Higher order equivalence of Bayes cross validation and 
WAIC. Springer Proceedings in Mathematics and Statistics, 
Information Geometry and Its Applications, pp.47-73, 2018. 

\bibitem{Watanabe2021}
Watanabe, S. WAIC and WBIC for mixture models. Behaviormetrika,
doi.org/10.1007/s41237-021-00133-z, 2021. 

\bibitem{Watanabe2021a}
Watanabe,S. 
Information criteria and cross validation for Bayesian inference in regular and singular cases. 
Japanese Journal of Statistics and Data Science volume,Vol 4, pp.1-19, 2021.

\bibitem{Watanabe2022}
Watanabe, S. Mathematical theory of Bayesian statistics where all models are wrong. Advancements in Bayesian Methods and Implementations, 
Handbook of statistics, Vol.47, pp.209-238, Elsevier, 2022.

\bibitem{Watanabe2022a}
Watanabe, S. Mathematical theory of Bayesian statistics
for unknown information source. to appear in 
Philosophical Transactions of the Royal Society A, arXiv:2206.05630, 2022. 

\bibitem{WatanabeT2022} 
Watanabe,T., Watanabe, S. 
Asymptotic behavior of Bayesian generalization error in multinomial mixtures.
arXiv:2203.06884. 

\bibitem{Wei2022}
Wei, S., Murfet, D., Gong, M., Li, H.,
Gell-Redman, J., Quella, T. 
Deep learning is singular, and That's good. 
IEEE Transactions on Neural Networks and Learning Systems, Vol.33,
pp.1-14, 2022. 

\bibitem{Yamazaki2003}
Yamazaki, K., Watanabe, S. Singularities in mixture models and upper bounds of stochastic complexity. International Journal of Neural Networks. 
Vol.16, No.7, pp.1029-1038, 2003.

\bibitem{Yamazaki2005}
Yamazaki, K.,  Watanabe,S. 
Algebraic geometry and stochastic complexity of hidden Markov models
Neurocomputing, Vol.69, pp.62-84, 2005. 

\bibitem{Yamazaki2005a}
Yamazaki, K.,  Watanabe, S. 
Singularities in complete bipartite graph-type boltzmann machines and upper bounds of stochastic complexities. 
IEEE transactions on neural networks, Vol.16, pp.312-324, 2005.

\bibitem{Yamazaki2007}
Yamazaki, K., Kawanabe, M., Watanabe, S., 
Sugiyama, M., M\"{u}ller, K.-R. 
Asymptotic bayesian generalization error when training and test distributions are different. 
Proceedings of the 24th international conference on Machine learning
pp. 1079-1086, 2007. 

\bibitem{Yamazaki2010}
Yamazaki, K., Aoyagi, M., Watanabe, S. 
Asymptotic analysis of Bayesian generalization error with Newton diagram. 
Neural Networks, Vol.23, pp.35-43, 2010. 

\bibitem{Yamazaki2016} 
Yamazaki, K. 
Asymptotic accuracy of Bayes estimation for latent variables with redundancy. Machine Learning . vol.102. pp.1-28, 2016.

\bibitem{Yamazaki2013}
Yamazaki, K., Kaji, D. 
Comparing two Bayes methods based on the free energy functions in Bernoulli mixtures.
Neural Networks. Vol.44, pp.36-43, 2013. 

\bibitem{Zwie2011}
Zwiernik, P.
An asymptotic behavior of the marginal likelihood for general Markov models.
The Journal of Machine Learning Research, vol.12, pp.3283-3310, 2011.


\end{thebibliography}
\end{document}